\documentclass[12pt]{article}
\usepackage{amsmath,amssymb,amscd,latexsym}

\numberwithin{equation}{section}
\usepackage[all]{xy}
\newcommand{\C}{{\mathbb{C}}}
\newcommand{\F}{\mathbb{F}}

\newcommand{\Q}{{\mathbb{Q}}}
\newcommand{\oQ}{\overline{\Q}}
\newcommand{\R}{{\mathbb{R}}}
\newcommand{\Z}{{\mathbb{Z}}}
\newcommand{\abb}{\mathrm{ab}}

\newcommand{\ddet}{\mathrm{det}}
\newcommand{\of}{\overline{f}}
\newcommand{\og}{\overline{g}}
\newcommand{\per}{\mathrm{per}\,}

\newcommand{\tr}{\mathrm{tr}}
\newcommand{\End}{\mathrm{End}\,}
\newcommand{\Fix}{\mathrm{Fix}\,}
\newcommand{\Gal}{\mathrm{Gal}\,}
\newcommand{\GL}{\mathrm{GL}}
\newcommand{\Ker}{\mathrm{Ker}\,}

\newcommand{\Nh}{{\mathcal N}}

\newcommand{\eb}{\mathfrak{b}}

\newcommand{\oA}{\overline{A}}

\newcommand{\tf}{\tilde{f}}

\newcommand{\tgamma}{\tilde{\gamma}}
\newcommand{\tGamma}{\tilde{\Gamma}}
\renewcommand{\mod}{\mathrm{mod}\,}

\newcommand{\silo}{\stackrel{\sim}{\longrightarrow}}
\newcommand{\tei}{\, | \,}

\newcommand{\halb}{\frac{1}{2}}
\newtheorem{theorem}{Theorem}[section]
\newtheorem{prop}[theorem]{Proposition}
\newtheorem{defn}[theorem]{Definition}
\newtheorem{cor}[theorem]{Corollary}
\newtheorem{example}[theorem]{Example}
\newtheorem{remark}[theorem]{Remark}
\newtheorem{claim}[theorem]{Claim}

\newenvironment{proofof}{\noindent {\it Proof of}}{\mbox{}\hspace*{\fill}$\Box$}
\newenvironment{proof}{\noindent {\bf Proof}}{\mbox{}\hspace*{\fill}$\Box$}
\begin{document}
\title{$p$-adic entropy and a $p$-adic Fuglede--Kadison determinant\\ 
{\normalsize \it Dedicated to Yuri Ivanovich Manin}}
\author{Christopher Deninger}

\maketitle

\section{Introduction} \label{sec:1}

In several instances the entropy $h (\varphi)$ of an automorphism $\varphi$ on a space $X$ can be calculated in terms of periodic points:
\begin{equation}\label{eq:1.1}
h (\varphi) = \lim_{n\to \infty} \frac{1}{n} \log |\Fix (\varphi^n)| \; .
\end{equation}
Here $\Fix (\varphi^n)$ is the set of fixed points of $\varphi^n$ on $X$. Let $\log_p : \Q^*_p \to \Z_p$ be the branch of the $p$-adic logarithm normalized by $\log_p (p) = 0$. The $p$-adic analogue of the limit \eqref{eq:1.1}, if it exists, may be viewed as a kind of entropy with values in the $p$-adic number field $\Q_p$,
\begin{equation} \label{eq:1.2}
h_p (\varphi) = \lim_{n \to \infty} \frac{1}{n} \log_p |\Fix (\varphi^n)| \; .
\end{equation}
It depends only on the action of $\varphi$ on $X$ viewed as a set. 

An earlier different approach to a $p$-adic entropy theory was mentioned to me by Amnon Besser. The usual definitions of measure theoretic or topological entropy have no obvious $p$-adic analogue since $\varlimsup$ or $\sup$ do not make sense $p$-adically and since the cardinalities of partitions, coverings and of separating or spanning sets do not behave reasonably in the $p$-adic metric.

Instead of actions of a single automorphism $\varphi$ we look more generally at actions of a countable discrete residually finite but not necessarily amenable group $\Gamma$ on a set $X$. Let us write $\Gamma_n \to e$ if $(\Gamma_n)$ is a sequence of cofinite normal subgroups of $\Gamma$ such that only the neutral element $e$ of $\Gamma$ is contained in infinitely many $\Gamma_n$'s. Let $\Fix_{\Gamma_n} (X)$ be the set of points in $X$ which are fixed by $\Gamma_n$. If the limit:
\begin{equation} \label{eq:1.3}
h_p := \lim_{n\to\infty} \frac{1}{(\Gamma : \Gamma_n)} \log_p |\Fix_{\Gamma_n} (X)|
\end{equation}
exists with respect to a choice of $\Gamma_n \to e$ we call it the $p$-adic entropy of the $\Gamma$-action on the set $X$ (with respect to the sequence $(\Gamma_n)$).

In this note we show that for an interesting class of $\Gamma$-actions the $p$-adic entropy exists independently of the choice of $\Gamma_n \to e$. In these examples $X$ is an abelian group and $\Gamma$ acts by automorphisms of groups. Namely, let $(\R / \Z)^{\Gamma}$ be the full shift on $\Gamma$ with values in the circle $\R / \Z$ and left $\Gamma$-action by $\gamma (x_{\gamma'}) = (x_{\gamma^{-1} \gamma'})$. For an element $f = \sum_{\gamma} a_{\gamma} \gamma$ in the integral group ring $\Z \Gamma$ consider the closed subshift $X_f \subset (\R / \Z)^{\Gamma}$ consisting of all sequences $(x_{\gamma'})$ which satisfy the equation
\[
\sum_{\gamma'} x_{\gamma'} a_{\gamma^{-1} \gamma'} = 0 \quad \mbox{in} \; (\R / \Z)^{\Gamma} \; \mbox{for all} \; \gamma \in \Gamma \; .
\]
In fact as in \cite{ER} we study more general systems defined by an $r \times r$-matrix over $\Z\Gamma$. However, in this introduction, for simplicity, we describe only the case $r = 1$. If $\Gamma$ is amenable, we denote by $h (f)$ the topological entropy of the $\Gamma$-action on $X_f$. 

The case $\Gamma = \Z^d$ is classical. Here we may view $f$ as a Laurent polynomial and according to \cite{LSW} the entropy is given by the (logarithmic) Mahler measure of $f$
\begin{equation} \label{eq:1.4}
h (f) = m (f) := \int_{T^d} \log |f (z)| \, d \mu (z) \; .
\end{equation}
Here $\mu$ is the normalized Haar measure on the $d$-torus $T^d$. According to \cite{LSW} the $\Z^d$-action on $X_f$ is expansive if and only if $f$ does not vanish in any point of $T^d$. By a theorem of Wiener this is also equivalent to $f$ being a unit in $L^1 (\Z^n)$. In this case $h (f)$ can be calculated in terms of periodic points, c.f. \cite{LSW} Theorem 7.1. See also \cite{S} for this theory.

What about a $p$-adic analogue? In \cite{D1} it was observed that in the expansive case $m (f)$ has an interpretation via the Deligne--Beilinson regulator map from algebraic $K$-theory to Deligne cohomology. Looking at the analogous regulator map from algebraic $K$-theory to syntomic cohomology one gets a suggestion what a (purely) $p$-adic Mahler measure $m_p (f)$ of $f$ should be, c.f. \cite{BD}. It can only be defined if $f$ does not vanish in any point of the $p$-adic $d$-torus $T^d_p = \{ z \in \C^d_p \tei |z_i|_p = 1 \} \; ,$ where $\C_p$ is the completion of a fixed algebraic closure $\oQ_p$ of $\Q_p$. In this case $m_p (f)$ is given by the convergent Snirelman integral
\begin{equation} \label{eq:1.5}
m_p (f) = \int_{T^d_p} \log_p f (z) \; .
\end{equation}
Recall that the Snirelman integral of a continuous function $F : T^d_p \to \C_p$ is defined by the following limit if it exists:
\[
\int_{T^d_p} F (z) := \lim_{N\to \infty \atop (N,p) = 1} \frac{1}{N^d} \sum_{\zeta \in \mu^d_N} F (\zeta) \; .
\]
Here $\mu_N$ is the group of $N$-th roots of unity in $\oQ^*_p$. 

For example, let $P (t) = a_m t^m + \ldots + a_r t^r$ be a polynomial in $\C_p [t]$ with $a_m , a_r \neq 0$ whose zeroes $\alpha$ satisfy $|\alpha|_p \neq 1$. Then, according to \cite{BD} Proposition 1.5 we have the following expression for the $p$-adic Mahler measure:
\begin{eqnarray} \label{eq:1.6}
m_p (f) & = & \log_p a_r - \sum_{0 < |\alpha|_p < 1} \log_p \alpha \\
 & = & \log_p a_m + \sum_{|\alpha|_p > 1} \log_p \alpha \; . \nonumber
\end{eqnarray}
For $d \ge 2$ there does not seem to be a simple formula for $m_p (f)$.

In \cite{BD} we mentioned the obvious problem to give an interpretation of $m_p (f)$ as a $p$-adically valued entropy. This is now provided by the following result:

\begin{theorem} \label{t:1.7}
Assume that $f \in \Z [\Z^d] = \Z [t^{\pm 1}_1 , \ldots , t^{\pm 1}_d]$ does not vanish in any point of the $p$-adic $d$-torus $T^d_p$. Then the $p$-adic entropy $h_p (f)$ of the $\Gamma = \Z^d$-action on $X_f$ in the sense of \eqref{eq:1.3} exists for all $\Gamma_n \to 0$ and we have $h_p (f) = m_p (f)$. 
\end{theorem}

Now we turn to more general groups $\Gamma$. In \cite{DS} extending \cite{D2} it was shown that for countable residually finite amenable groups $\Gamma$ and elements $f$ in $\Z \Gamma$ which are invertible in $L^1 (\Gamma)$ we have
\begin{equation} \label{eq:1.8}
h (f) = \log \ddet_{\Nh\Gamma} f \; .
\end{equation}
Here $\det_{\Nh\Gamma}$ is the Fuglede--Kadison determinant \cite{FK} on the units of the von~Neumann algebra $\Nh\Gamma \supset L^1 \Gamma \supset \Z\Gamma$ of $\Gamma$. In fact, equation \eqref{eq:1.8} holds without the condition of amenability if $h (f)$ is replaced by the quantity:
\[
h_{\per} (f) := \lim_{n\to \infty} \frac{1}{(\Gamma : \Gamma_n)} \log |\Fix_{\Gamma_n} (X)| \; .
\]
For the $\Gamma$-action on $X_f$ this limit exists and is independent of the choice of sequence $\Gamma_n \to e$.

In the $p$-adic case, instead of working with a $p$-adic $L^1$-convolution algebra it is more natural to work with the bigger convolution algebra $c_0 (\Gamma)$. It consists of all formal series $x = \sum_{\gamma} x_{\gamma} \gamma$ with $x_{\gamma} \in \Q_p$ and $|x_{\gamma}|_p \to 0$ as $\gamma \to \infty$ in $\Gamma$.

For $\Gamma = \Z^d$ it is known that $f \in \Z [\Z^d]$ does not vanish in any point of the $p$-adic $d$-torus $T^d_p$ if and only if $f$ is a unit in the algebra $c_0 (\Z^d)$. Hence in general, it is natural to look for a $p$-adic analogue of formula \eqref{eq:1.8} for all $f \in \Z \Gamma$ which are units in $c_0 (\Gamma)$. In the $p$-adic case there is no analogue for the theory of von~Neumann algebras and for the functional calculus used to define $\det_{\Nh\Gamma}$. However using some algebraic $K$-theory and the results of \cite{FL}, \cite{L} and \cite{KLM} we can define a $p$-adic analogue $\log_p \det_{\Gamma}$ of $\log \det_{\Nh\Gamma}$ for suitable classes of groups $\Gamma$. For example we get the following result generalizing theorem \ref{t:1.7}:

\begin{theorem} \label{t:1.9}
Assume that the residually finite group $\Gamma$ is elementary amenable and torsion-free. Let $f$ be an element of $\Z\Gamma$ which is a unit in $c_0 (\Gamma)$. Then the $p$-adic entropy $h_p (f)$ of the $\Gamma$-action on $X_f$ in the sense of \eqref{eq:1.3} exists for all $\Gamma_n \to e$ and we have
\[
h_p (f) = \log_p \ddet_{\Gamma} f \; .
\]
\end{theorem}

Acknowledgement: I would like to thank my colleagues Wolfgang L\"uck and Peter Schneider for helpful conversations.
\section{Preliminaries} \label{sec:2}

Fix an integer $r \ge 1$ and set $T = T^r = (\R / \Z)^r$. For a discrete group $\Gamma$ let $T^{\Gamma}$ be the full shift with left $\Gamma$-action by $\gamma (x_{\gamma'}) = (x_{\gamma^{-1} \gamma'})$. Write $M_r (R)$ for the ring of $r \times r$-matrices over a ring $R$. For an element $f = \sum a_{\gamma} \gamma$ in $M_r (\Z) [\Gamma] = M_r (\Z\Gamma)$ the closed subshift $X_f \subset T^{\Gamma}$ is defined as the closed subgroup consisting of all sequences with
\[
\sum_{\gamma'} x_{\gamma'} a^*_{\gamma^{-1} \gamma'} = 0 \quad \mbox{in} \; T^{\Gamma} \; \mbox{for all} \; \gamma \in \Gamma \; .
\]
Here $a^*$ denotes the transpose of a matrix $a$ in $M_r (\Z)$. The group ring $M_r (\Z) [\Gamma]$ is equipped with an anti-involution $*$ defined by $f^* = \sum_{\gamma} a^*_{\gamma^{-1}} \gamma$ for $f = \sum_{\gamma} a_{\gamma} \gamma$. 

Let $\rho_f$ be right multiplication by $f^*$ on the group $T [[\Gamma]]$ of formal $T$-valued series on $\Gamma$. For $x = \sum_{\gamma} x_{\gamma} \gamma$ in $T [[\Gamma]]$ we have
\[
\rho_f (x) = \sum_{\gamma} x_{\gamma} \gamma \sum_{\gamma} a^*_{\gamma^{-1}} \gamma = \sum_{\gamma} \Big( \sum_{\gamma'} x_{\gamma'} a^*_{\gamma^{-1} \gamma'} \Big) \gamma \; .
\]
Hence we see that 
\[
X_f = \Ker (\rho_f : T [[\Gamma]] \longrightarrow T [[\Gamma]]) \; 
\]
where on the right hand side the group $\Gamma$ acts by left multiplication. Let $N$ be a normal subgroup of $\Gamma$ with quotient map $\sim : \Gamma \to \tilde{\Gamma} = \Gamma / N$. Set
\[
\tilde{f} = \sum_{\gamma} a_{\gamma} \tilde{\gamma} = \sum_{\delta \in \tilde{\Gamma}} \Big( \sum_{\gamma \in \delta} a_{\gamma} \Big) \delta \quad \mbox{in} \; M_r (\Z) [\tilde{\Gamma}] \; .
\]
This is the image of $f$ under the reduction map $M_r (\Z) [\Gamma] \to M_r (\Z) [\tilde{\Gamma}]$. Under the natural isomorphism
\[
T [[\tilde{\Gamma}]] \silo \Fix_N (T [[\Gamma]])
\]
mapping $\sum_{\delta} x_{\delta} \delta$ to $\sum_{\gamma} x_{\tilde{\gamma}} \gamma$ the action $\rho_{\tilde{f}}$ corresponds to the restriction of $\rho_f$. Hence we have
\[
\Fix_N (X_f) = \Ker (\rho_{\tilde{f}} : T [[\tilde{\Gamma}]] \longrightarrow T [[\tilde{\Gamma}]]) = X_{\tf} \; .
\]
If we assume that $\tilde{\Gamma}$ is finite we get that
\[
\Fix_N (X_f) = \rho^{-1}_{\tilde{f},\R} (\Z \tilde{\Gamma})^r / (\Z \tilde{\Gamma})^r
\]
for the endomorphism $\rho_{\tilde{f},\R}$ of right multiplication by $\tilde{f}^*$ on $(\R \tilde{\Gamma})^r$. This implies the following fact, c.f. \cite{DS}, Corollary 4.3:

\begin{prop} \label{t:2.1}
Let $\tilde{\Gamma}$ be finite. Then $\rho_{\tilde{f}}$ is an isomorphism of $(\Q \tilde{\Gamma})^r$ if and only if $\Fix_N (X_f)$ is finite. In this case the order is given by
\[
|\Fix_N (X_f)| = \pm \ddet \rho_{\tilde{f}} \; .
\]
\end{prop}

This follows from the fact that for an isomorphism $\varphi$ of a finite dimensional real vector space $V$ and a lattice $\Lambda$ in $V$ with $\varphi (\Lambda) \subset \Lambda$ we have:
\[
|\varphi^{-1} \Lambda / \Lambda| = |\Lambda / \varphi (\Lambda)| = |\det (\varphi \tei V)|\; .
\]

For any countable discrete group $\Gamma$ let $c_0 (\Gamma)$ be the set of formal series $\sum_{\gamma} x_{\gamma} \gamma$ with $x_{\gamma} \in \Q_p$ and $|x_{\gamma}|_p \to 0$ for $\gamma \to \infty$. This means that for any $\varepsilon > 0$ there is a finite subset $S \subset \Gamma$ such that $|x_{\gamma}|_p < \varepsilon$ for all $\gamma \in \Gamma \setminus S$. The set $c_0 (\Gamma)$ is a $\Q_p$-vector space and it becomes a $\Q_p$-algebra with the product
\begin{equation} \label{eq:2.2}
\sum_{\gamma} x_{\gamma} \gamma \cdot \sum_{\gamma} y_{\gamma} \gamma = \sum_{\gamma} \Big( \sum_{\gamma' \gamma'' = \gamma} x_{\gamma'} y_{\gamma''} \Big) \gamma \; .
\end{equation}
Note that the sums
\[
\sum_{\gamma' \gamma'' = \gamma} x_{\gamma'} y_{\gamma''} = \sum_{\gamma'} x_{\gamma'} y_{\gamma^{'-1} \gamma}
\]
converge $p$-adically for every $\gamma$ since $\lim_{\gamma' \to \infty} |x_{\gamma'} y_{\gamma^{'-1} \gamma}|_p = 0$. The value is independent of the order of summation. Moreover, because of the inequality
\begin{equation} \label{eq:2.3}
\Big| \sum_{\gamma'} x_{\gamma'} y_{\gamma^{'-1} \gamma} \Big|_p \le \sup_{\gamma'} |x_{\gamma'} y_{\gamma^{'-1} \gamma}|_p \; ,
\end{equation}
we have
\[
\lim_{\gamma \to \infty} \sum_{\gamma'} x_{\gamma'} y_{\gamma^{'-1} \gamma} = 0 \; ,
\]
so that the product \eqref{eq:2.2} is well defined. We may also view $c_0 (\Gamma)$ as an algebra of $\Q_p$-valued functions on $\Gamma$ under convolution. 

The $\Q_p$-algebra $c_0 (\Gamma)$ is complete in the norm
\[
\| \sum_{\gamma} x_{\gamma} \gamma \| = \sup_{\gamma} |x_{\gamma}|_p = \max_{\gamma} |x_{\gamma} |_p \; .
\]
The norm satisfies the following properties:
\begin{align}
 & \| x \| = 0 \quad \mbox{if and only if} \; x = 0 \label{eq:2.4} \\
 & \| x + y \| \le \max (\|x \| , \| y \|) \label{eq:2.5} \\
& \| \lambda x \| = |\lambda|_p \, \| x \| \quad \mbox{for all} \; \lambda \in \Q_p \label{eq:2.6} \\
& \|xy \| \le \| x \| \, \| y \|  \quad \mbox{and} \; \| 1 \| = 1 \label{eq:2.7}
\end{align}
Hence $c_0 (\Gamma)$ is a $p$-adic Banach algebra over $\Q_p$, i.e. a unital $\Q_p$-algebra $B$ which is complete with respect to a norm $\|\ \| : B \to \R^{\ge 0}$ satisfying conditions \eqref{eq:2.4}--\eqref{eq:2.7}. 

We will only consider Banach algebras where $\| \; \|$ takes values in $p^{\Z} \cup \{ 0 \}$. The subring $A = B^0$ of elements $x$ in $B$ of norm $\| x \| \le 1$ is a $p$-adic Banach algebra over $\Z_p$, defined similarly as before. An example is given by
\[
c_0 (\Gamma , \Z_p) = c_0 (\Gamma)^0 = \{ \sum x_{\gamma} \gamma \tei x_{\gamma} \in \Z_p \; \mbox{with} \; \lim_{\gamma \to \infty} |x_{\gamma}|_p = 0 \} \; .
\]
In this case the residue algebra $A / pA$ over $\F_p$ is isomorphic to the group ring of $\Gamma$ over \nolinebreak $\F_p$:
\begin{equation} \label{eq:2.8}
c_0 (\Gamma , \Z_p) / p c_0 (\Gamma , \Z_p) = \F_p [\Gamma] \; .
\end{equation}
The $1$-units $U^1 = 1 + pA$ form a subgroup of $A^*$ since
\[
(1 + pa)^{-1} := \sum^{\infty}_{\nu =0} (-pa)^{\nu}
\]
provides an inverse of $1 + pa \in U^1$ in $U^1$. It is easy to see that one has an exact sequence of groups
\begin{equation} \label{eq:2.9}
1 \longrightarrow U^1 \longrightarrow A^* \longrightarrow (A / pA)^* \longrightarrow 1 \; .
\end{equation}
For $A = c_0 (\Gamma , \Z_p)$ this is the exact seqence
\begin{equation} \label{eq:2.10}
1 \longrightarrow 1 + pc_0 (\Gamma , \Z_p) \longrightarrow c_0 (\Gamma , \Z_p)^* \longrightarrow \F_p [\Gamma]^* \longrightarrow 1 \; .
\end{equation}
Concerning the units of a $p$-adic Banach algebra over $\Q_p$, we have the following known fact:

\begin{prop} \label{t:2.11}
Let $B$ be a $p$-adic Banach algebra over $\Q_p$ whose norm takes values in $p^{\Z} \cup \{ 0 \}$ and set $A = B^0$. If the residue algebra $A / pA$ has no zero divisors, then we have
\[
B^* = p^{\Z} A^* \quad \mbox{and} \quad p^{\Z} \cap A^* = 1 \; .
\]
\end{prop}

\begin{proof}
For $f$ in $B^*$ set $g = 1/f$. Let $\nu , \mu$ be such that $f_1 = p^{\nu} f$ and $g_1 = p^{\mu} g$ have norm one. The reductions $\of_1 , \og_1$ of $f_1 , g_1$ are non-zero. In the equation $f_1 g_1 = p^{\nu + \mu}$ we have $\nu + \mu \ge 0$. Reducing $\mod p$ we find that $0 \neq \of_1  \og_1 = p^{\nu + \mu} \, \mod p$ in $A / pA$. Hence we have $\nu + \mu = 0$ and therefore $f_1 g_1 = 1$. The first assertion follows. Because of \eqref{eq:2.7} we have $\| a \| = 1$ for $a \in A^*$ and $\| p^{\nu} \| = p^{-\nu}$. This implies the second assertion.
\end{proof}

\begin{example} \label{t:2.12}
\em For the group $\Gamma = \Z^d$ the algebra $c_0 (\Z^d)$ can be identified with the affinoid commutative algebra $\Q_p \langle t^{\pm 1}_1 , \ldots , t^{\pm 1}_d \rangle$ of power series $\sum_{\nu \in \Z^d} x_{\nu} t^{\nu}$ with $x_{\nu} \in \Q_p$ and $\lim_{|\nu| \to \infty} |x_{\nu}|_p = 0$. Note that these power series can be viewed as functions on $T^d_p$. The residue algebra is $\F_p [\Z^d] = \F_p [t^{\pm 1}_1 , \ldots , t^{\pm 1}_d]$. It has no zero divisors and its groups of units is
\[
\F_p [\Z^d]^* = \F^*_p t^{\Z}_1 \cdots t^{\Z}_d \; .
\]
The preceeding proposition and the exact sequence \eqref{eq:2.10} now give a decomposition into a direct product of groups
\[
c_0 (\Z^d)^* = p^{\Z} \mu_{p-1} t^{\Z}_1 \cdots t^{\Z}_d (1 + p\,c_0 (\Z^d , \Z_p)) \; .
\]
\end{example}

\begin{prop} \label{t:2.13}
For $f$ in $\Q_p [\Z^d] = \Q_p [t^{\pm 1}_1 , \ldots , t^{\pm 1}_d]$ the following properties are equivalent:\\
a) We have $f (z) \neq 0$ for every $z$ in $T^d_p$\\
b) $f$ is a unit in $c_0 (\Z^d)^*$\\
c) $f$ has the form $f (t) = c t^{\nu} (1 + pg (t))$ for some $c \in \Q^*_p$, $\nu \in \Z^d$ and $g (t)$ in $c_0 (\Z^d , \Z_p)$.
\end{prop}

\begin{proof}
We have seen that b) and c) are equivalent and it is clear that both b) and c) imply a). For proving that a) implies b) note that the maximal ideals of $c_0 (\Z^d) = \Q_p \langle t^{\pm 1}_1 , \ldots , t^{\pm 1}_d \rangle$ correspond to the orbits of the $\Gal (\oQ_p / \Q_p)$-operation on $T^d_p \cap (\oQ^*_p)^d$. Hence $f$ is not contained in any maximal ideal of $c_0 (\Z^d)$ by assumption a) and therefore $f$ is a unit.
\end{proof}

\section{The Frobenius group determinant and a proof of theorem \ref{t:1.7}} \label{sec:3}

The map $\Z \Gamma \to \Z \tGamma$ for $\tGamma = \Gamma / N$ from the beginning of the last section can be extended to a homomorphism of $\Q_p$-algebras $c_0 (\Gamma) \to c_0 (\tGamma)$ by sending $f = \sum a_{\gamma} \gamma$ to $\tf = \sum a_{\gamma} \tgamma$. Note that this is well defined by the ultrametric inequality and that we have $\| \tf \|\le \| f \|$. If $\tGamma$ is finite, we have $c_0 (\tGamma) = \Q_p \tGamma$ and hence we obtain a homomorphism of groups $\GL_r (c_0 (\Gamma)) \to \GL_r (\Q_p\tGamma)$. It follows that for $f$ in $M_r (\Z\Gamma) \cap \GL_r (c_0 (\Gamma))$ the endomorphism $\rho_{\tf}$ of $(\Q \tGamma)^r$ is an isomorphism. Together with proposition \ref{t:2.1} we have shown the first equation in the following proposition: 
\begin{prop} \label{t:3.1}
Let $\Gamma$ be a discrete group and $N$ a normal subgroup with finite quotient group $\tGamma$. For $f$ in $M_r (\Z \Gamma) \cap \GL_r (c_0 (\Gamma))$ the set $\Fix_N (X_f)$ is finite and we have
\begin{eqnarray*}
|\Fix_N (X_f)| & = & \pm \det \rho_{\tf} \\
 & = & \pm \prod_{\pi} \ddet_{\oQ_p} \Big( \sum_{\gamma} a^*_{\gamma} \otimes \rho_{\pi} (\tgamma) \Big)^{d_{\pi}} \; .
\end{eqnarray*}
Here $\pi$ runs over the equivalence classes of irreducible representations $\rho_{\pi}$ of $\tGamma$ on $\oQ_p$-vector spaces $V_{\pi}$ and $d_{\pi}$ is the degree $\dim V_{\pi}$ of $\pi$.
\end{prop}

\begin{proof}
It remains to prove the second equation which is essentially due to Frobenius. Consider $\oQ_p \tGamma$ as a representation of $\tGamma$ via the map $\tgamma \mapsto \rho (\tgamma) :=$ right multiplication with $\tgamma^{-1}$. This (``right regular'') representation decomposes as follows into irreducible representations c.f. \cite{S} I, 2.2.4
\[
\oQ_p \tGamma \cong \bigoplus_{\pi} V^{d_{\pi}}_{\pi} \; .
\]
The endomorphism
\[
\rho_{\tf} = \sum_{\gamma} a^*_{\gamma} \otimes \rho (\tgamma)
\]
on 
\[
(\oQ_p \tGamma)^r = \oQ^r_p \otimes \oQ_p \tGamma
\]
therefore corresponds to the endomorphism
\[
\bigoplus_{\pi} \Big( \sum_{\gamma} a^*_{\gamma} \otimes \rho_{\pi} (\tgamma)\Big)^{d_{\pi}} \quad \mbox{on} \quad \bigoplus_{\pi} \oQ^r_p \otimes V^{d_{\pi}}_{\pi} = \bigoplus_{\pi} (\oQ^r_p \otimes V_{\pi})^{d_{\pi}} \; .
\]
Hence the formula follows.
\end{proof}

{\bf Remark} In the real case and for the Heisenberg group, Klaus Schmidt previously used the group determinant to calculate $|\Fix_{\Gamma_n} (X_f)|$ for $f$ in $L^1 (\Gamma)^*$ and certain $\Gamma_n$.

The following result generalizes theorem \ref{t:1.7} from the introduction, at least for a particular sequence $\Gamma_n \to 0$:

\begin{theorem} \label{t:3.2}
Let $f = \sum_{\nu \in \Z^n} a_{\nu} t^{\nu}$ in $M_r (\Z [t^{\pm 1}_1 , \ldots , t^{\pm 1}_d])$ be invertible in every point of the $p$-adic $d$-torus $T^d_p$. Then the $p$-adic entropy $h_p (f)$ of the $\Gamma = \Z^d$-action on $X_f$ exists in the sense of \eqref{eq:1.3} for the sequence $\Gamma_n = (n \Z)^d \to 0$ with $n$ prime to $p$, and we have
\[
h_p (f) = m_p (\det f) \; .
\]
\end{theorem}

\begin{proof}
By assumption, the Laurent polynomial $\det f$ does not vanish in any point of $T^d_p$. Hence $\det f$ is a unit in $c_0 (\Gamma)$ by proposition \ref{t:2.13}. It follows from proposition \ref{t:3.1} that we have:
\[
|\Fix_{\Gamma_n} (X_f)| = \pm \prod_{\chi} \ddet_{\oQ_p} \Big( \sum_{\nu \in \Z^d} a^*_{\nu} \otimes \chi (\nu) \Big) 
\]
where $\chi$ runs over the characters of $\Gamma / \Gamma_n = (\Z / n \Z)^d$. These correspond via $\chi (\nu) = \zeta^{\nu}$ to the elements $\zeta$ of $\mu^d_n$. Viewing $f$ as a matrix of functions on $T^d_p$ we therefore get the formulas
\begin{eqnarray*}
|\Fix_{\Gamma_n} (X_f)| & = & \pm \prod_{\zeta \in \mu^d_n} \ddet_{\oQ_p} \Big( \sum_{\nu \in \Z^d} a^*_{\nu} \zeta^{\nu} \Big) \\
 & = & \pm \prod_{\zeta \in \mu^d_n} \ddet_{\oQ_p} (f (\zeta)) \\
 & = & \pm \prod_{\zeta \in \mu^d_n} (\det f) (\zeta) \; .
\end{eqnarray*}
Thus the $p$-adic entropy $h_p (f)$ of the $\Gamma$-action on $X_f$ with respect to the above sequence is given by:
\begin{eqnarray*}
h_p (f) & = & \lim_{n \to \infty \atop (n,p) = 1} \frac{1}{(\Gamma : \Gamma_n)} \log_p |\Fix_{\Gamma_n} (X_f)| \\
& = & \lim_{n\to \infty \atop (n,p) = 1} \frac{1}{n^d} \sum_{\zeta \in \mu^d_n} \log_p (\det f) (\zeta) \\
& = & \int_{T^d_p} \log_p \det f = m_p (\det f) \; .
\end{eqnarray*}
Note here that for the Laurent polynomials $\det f$ under consideration the Snirelman integral exists by \cite{BD} Proposition 1.3.
\end{proof}

{\bf Remark} A suitable generalization of that proposition would give theorems \ref{t:3.2} and \ref{t:1.7} for general sequences $\Gamma_n \to 0$ in $\Gamma = \Z^d$. We leave this to the interested reader since the general case of theorem \ref{t:1.7} is also a corollary of theorem \ref{t:1.9} which will be proved by a different method in section \ref{sec:5}.

{\bf Example} The polynomial in one variable $f (T) = 2T^2 - T + 2$ does not vanish in any point of the $2$-adic circle $T^1_2$. In this sense, $X_f$ is ``$2$-adically expansive''. Consider the square root of $-15$ in $\Z_2$ given by the $2$-adically convergent series
\[
\sqrt{-15} = (1 + (-16))^{1/2} = \sum^{\infty}_{\nu = 0} {1/2 \choose \nu} (-1)^{\nu} 2^{4\nu} \; .
\]
The zeroes of $f (T)$ in $\oQ_2$ are given by $\alpha_{\pm} = \frac{1}{4} (1 \pm \sqrt{-15}) \in \Q_2$ where $|\alpha_+|_2 = 2$ and $|\alpha_-|_2 = 1/2$. Successive approximations for $\alpha_+$ coming from the series for $\sqrt{-15}$ are $1/2 , -3/2 , -19/2 , -83/2$. By theorem \ref{t:3.2} and formula \eqref{eq:1.6} the $2$-adic entropy of $X_f$ is given by
\[
h_2 (f) = \log_2 \alpha_+ \in \Z_2 \; .
\]
Note that $f$ viewed as a complex valued function has both its zeroes on $S^1$, so that $X_f$ is not expansive in the usual sense. The topological entropy is $h (f) = \log 2$. 

\section{The logarithm on the $1$-units of a $p$-adic Banach algebra} \label{sec:4}

For a discrete group $\Gamma$ we would like to define a homomorphism
\[
\log_p \ddet_{\Gamma} : c_0 (\Gamma)^* \longrightarrow \Q_p
\]
which should be a $p$-adic replacement for the map
\[
\log \ddet_{\Nh\Gamma} : L^1 (\Gamma)^* \subset (\Nh\Gamma)^* \longrightarrow \R \; .
\]
More generally, we would like to define such a map on $\GL_r (c_0 (\Gamma))$. In this section we give the definition on the subgroup of $1$-units and relate $\log_p \det_{\Gamma}$ to $p$-adic entropy. The extension to a map on all of $c_0 (\Gamma)^*$ will be done in the next section for suitable classes of groups $\Gamma$ using rather deep facts about group rings.

Let $B$ be a $p$-adic Banach algebra over $\Q_p$ whose norm $\| \; \|$ takes values in $p^{\Z} \cup \{ 0 \}$. A trace functional on $B$ is a continuous linear map $\tr_B : B \to \Q_p$ which vanishes on commutators $[a,b] = ab - ba$ of elements in $B$. For $b \in B$ and $c \in B^*$ we have
\begin{equation} \label{eq:4.1}
\tr_B (cbc^{-1}) = \tr_B (b) \; .
\end{equation}
Set $A = B^0 = \{ b \in B \tei \| b \| \le 1 \}$ and let $U^1$ be the normal subgroup of $1$-units in $A^*$. The logarithmic series
\[
\log : U^1 \longrightarrow A \; , \; \log u = - \sum^{\infty}_{\nu = 1} \frac{(1 - u)^{\nu}}{\nu}
\]
converges and defines a continuous map. An argument with formal power series shows that we have
\begin{equation} \label{eq:4.2}
\log uv = \log u + \log v 
\end{equation}
if the elements $u$ and $v$ in $U^1$ {\it commute with each other}. 

The next result is a consequence of the Campbell--Baker--Hausdorff formula. 

\begin{theorem} \label{t:4.3}
The map
\[
\tr_B \log : U^1 \longrightarrow \Z_p
\]
is a homomorphism. For $u$ in $U^1$ and $a$ in $A^*$ we have
\begin{equation} \label{eq:4.4}
\tr_B \log (a u a^{-1}) = \tr_B \log (u) \; .
\end{equation}
\end{theorem}

\begin{proof}
Formula \eqref{eq:4.4} follows from \eqref{eq:4.1}. From \cite{B} Ch. II, \S\,8 we get the following information about $\log$. Set $G = \{ b \in B \tei \| b \| < p^{-\frac{1}{p-1}} \}$. Then the exponential series defines a bijection $\exp : G \silo 1+G$ with inverse $\exp^{-1} = \log \, |_{1+G}$. For $x,y$ in $G$ we have
\begin{equation} \label{eq:4.5}
\exp x \cdot \exp y = \exp h (x,y)
\end{equation}
where $h (x,y) \in G$ is given by a convergent series in $B$. It has the form
\[
h (x,y) = x + y + \mbox{series of (iterated) commutators} \; .
\]
Elements $u,v$ of $1+G$ have the form $u = \exp x$ and $v = \exp y$. Taking the $\log$ of relation \eqref{eq:4.5} and applying $\tr_B$ we get
\begin{eqnarray}
\tr_B \log (uv) & = & \tr_B h (x,y) = \tr_B (x+y) \nonumber \\
 & = & \tr_B \log u + \tr_B \log v \; . \label{eq:4.6}
\end{eqnarray}
Hence $\tr_B \log$ is a homomorphism on the subgroup $1+G$ of $U^1$. By assumption the norm of $B$ takes values in $p^{\Z} \cup \{ 0 \}$. For $p \neq 2$ we therefore have $1 + G = U^1$ and we are done. \\
For $p = 2$ the restriction of the map
\[
\varphi = \tr_B \log : U^1 \to \Q_2
\]
to $1 +G = 1 + 4A$ is a homomorphism by \eqref{eq:4.6}. We have to show that it is a homomorphism on $U^1 = 1 + 2A$ as well. For $u$ in $U^1$ we have $\varphi (u) = \halb \varphi (u^2)$ by \eqref{eq:4.2} and $u^2$ lies in $1 + 4A$. Now consider elements $u, v$ in $U^1$. Then we have
\begin{eqnarray*}
\varphi (uv) & = & \halb \varphi ((uv)^2) = \halb \varphi (uvuv) \overset{\eqref{eq:4.4}}{=} \halb \varphi (u^2 v u v u^{-1}) \\
 & = & \halb \varphi (u^2) + \halb \varphi (vuvu^{-1})
\end{eqnarray*}
since $u^2$ and $v uvu^{-1}$ lie in $1 + 4A$ where $\varphi$ is a homomorphism. By similar arguments we get
\begin{eqnarray*}
\varphi (uv) & = & \varphi (u) + \halb \varphi (v^2 uvu^{-1} v^{-1}) \\
& = & \varphi (u) + \varphi (v) + \halb \varphi (uvu^{-1} v^{-1}) \; .
\end{eqnarray*}
Thus we must show that $\varphi (uvu^{-1} v^{-1}) = 0$. By \eqref{eq:4.4} we have $\varphi (uvu^{-1} v^{-1}) = \varphi (vu^{-1} v^{-1} u)$ and hence using that both $uvu^{-1} v^{-1}$ and $v u^{-1} v^{-1} u$ lie in $1 + 4A$ we find
\begin{eqnarray*}
2 \varphi (uv u^{-1} v^{-1}) & = & \varphi (uvu^{-1} v^{-1}) + \varphi (vu^{-1} v^{-1} u) \\
 & = & \varphi (uv u^{-2} v^{-1} u) \overset{\eqref{eq:4.4}}{=} \varphi (u^{-2} v^{-1} u^2 v) \\
& = & \varphi (u^{-2}) + \varphi (v^{-1} u^2 v) \\
& \overset{\eqref{eq:4.4}}{=} & \varphi (u^{-2}) + \varphi (u^2) = \varphi (e) \\
 & = & 0 \; .
\end{eqnarray*}
\end{proof}

For a discrete group $\Gamma$, the map
\[
\tr_{\Gamma} : c_0 (\Gamma) \longrightarrow \Q_p \; , \; \tr_{\Gamma} ({\textstyle \sum} a_{\gamma} \gamma) = a_e
\]
defines a trace functional on $c_0 (\Gamma)$. Let $B = M_r (c_0 (\Gamma))$ be the $p$-adic Banach algebra over $\Q_p$ of $r \times r$-matrices $(a_{ij})$ with entries in $c_0 (\Gamma)$ and equipped with the norm $\| (a_{ij}) \| = \max_{ij} \| a_{ij} \|$. The composition:
\[
\tr_{\Gamma} : M_r (c_0 (\Gamma)) \xrightarrow{\tr} c_0 (\Gamma) \xrightarrow{\tr_{\Gamma}} \Q_p
\]
defines a trace functional on $M_r (c_0 (\Gamma))$. 

The algebra $A = B^0$ is given by $M_r (c_0 (\Gamma , \Z_p))$ and we have $U^1 = 1 + pM_r (c_0 (\Gamma , \Z_p))$. The exact sequence \eqref{eq:2.9} becomes the exact sequence of groups:
\begin{equation} \label{eq:4.7}
1 \longrightarrow 1 + p M_r (c_0 (\Gamma , \Z_p)) \longrightarrow \GL_r (c_0 (\Gamma , \Z_p)) \longrightarrow \GL_r (\F_p \Gamma) \longrightarrow 1 \; .
\end{equation}
According to theorem \ref{t:4.3} the map
\begin{equation} \label{eq:4.8}
\log_p \ddet_{\Gamma} := \tr_{\Gamma} \log : 1 + p M_r (c_0 (\Gamma , \Z_p)) \longrightarrow \Z_p
\end{equation}
is a homomorphism of groups.

\begin{example} \label{t:4.9}
\rm For $\Gamma = \Z^d$, in the notation of example \ref{t:2.12} we have a commutative diagram, c.f. \cite{BD} Lemma 1.1
\[
\xymatrix{
c_0 (\Gamma) \ar[r]^-{\sim} \ar[d]_{\tr_{\Gamma}}&  \Q_p \langle t^{\pm 1}_1 , \ldots , t^{\pm 1}_d \rangle \ar[d]^{\int_{T^d_p}} \\
\Q_p \ar@{=}[r] & \Q_p \; .
}
\]
It follows that for a $1$-unit $f$ in $M_r (c_0 (\Gamma))$ we have:
\begin{equation} \label{eq:4.10}
\log_p \ddet_{\Gamma} f = \int_{T^d_p} \log \det f = m_p (\det f) \; .
\end{equation}
Here we have used the relation
\begin{equation} \label{eq:4.11}
\tr \log f = \log \det f \quad \mbox{in} \; c_0 (\Gamma) \; ,
\end{equation}
where $\det : \GL_r (c_0 (\Gamma)) \to c_0 (\Gamma)^*$ is the determinant and $\tr$ the trace for matrices over the {\it commutative} ring $c_0 (\Gamma)$. Note that $\det$ maps $1$-units to $1$-units. Relation \eqref{eq:4.11} can be proved by embedding the integral domain $c_0 (\Gamma) = \Q_p \langle t^{\pm 1}_1 , \ldots , t^{\pm 1}_d \rangle$ into its quotient field and applying \cite{H} Appendix C, Lemma 4.1.
\end{example}

For finite groups $\Gamma$ the map $\log_p \det_{\Gamma}$ can be calculated as follows. For $f$ in $M_r (c_0 (\Gamma)) = M_r (\Q_p \Gamma)$ let $\rho_f$ be the endomorphism of $(\Q_p \Gamma)^r$ by right multiplication with $f^*$ and $\det_{\Q_p} (\rho_f)$ its determinant over $\Q_p$. 

\begin{prop} \label{t:4.12}
Let $\Gamma$ be finite. Then we have
\begin{equation} \label{eq:4.13}
\log_p \ddet_{\Gamma} f = \frac{1}{|\Gamma|} \log_p \ddet_{\Q_p} (\rho_f)
\end{equation}
for $f$ in $1 + p M_r (\Z_p \Gamma)$.
\end{prop}

{\bf Remark} Since $\rho_{fg} = \rho_f \rho_g$, the equation in the proposition shows that $\log_p \det_{\Gamma}$ is a homomorphism -- something we know in general by theorem \ref{t:4.3}. For finite $\Gamma$ the group $\GL_r (\F_p \Gamma)$ is finite. Hence, by \eqref{eq:4.7} there is at most one way to extend $\log_p \det_{\Gamma}$ from $1 + p M_r (\Z_p \Gamma)$ to a homomorphism from $\GL_r (\Z_p \Gamma)$ to $\Q_p$. Namely, we have to set
\[
\log_p \ddet_{\Gamma} f := \frac{1}{N} \log_p \ddet_{\Gamma} f^N \; ,
\]
where $N \ge 1$ is any integer with $\of^N = 1$ in $\GL_r (\F_p \Gamma)$. Because of \eqref{eq:4.2} this is well defined but it is not clear from the definition that we get a homomorphism. However, for the same $f , N$ we have
\[
\log_p \ddet_{\Q_p} (\rho_f) = \frac{1}{N} \log_p \ddet_{\Q_p} (\rho_{f^N}) \; .
\]
Hence equation \eqref{eq:4.13} holds for all $f$ in $\GL_r (\Z_p \Gamma)$ and it follows that $\log_p \ddet_{\Gamma}$ extends to a {\it homomorphism} on $\GL_r (\Z_p \Gamma)$. In the next section such arguments will be generalized to infinite groups with the help of $K$-theory.

\begin{proofof} {\bf \ref{t:4.12}}
Under the continuous homomorphism of $p$-adic Banach algebras over $\Z_p$
\[
\rho : M_r (\Z_p \Gamma) \longrightarrow \End_{\Z_p} (\Z_p \Gamma)^r
\]
the groups of $1$-units are mapped to each other. Hence we have
\begin{equation} \label{eq:4.14}
\log \rho_f = \rho_{\log f}
\end{equation}
for $f$ in $1 + pM_r (\Z_p \Gamma)$. 

On the other hand we have
\begin{equation} \label{eq:4.15}
\tr_{\Gamma} (g) = \frac{1}{|\Gamma|} \tr (\rho_g)
\end{equation}
for any element $g$ of $M_r (\Q_p \Gamma)$. This is proved first for $r = 1$ by checking the cases where $g = \gamma$ is an element of $\Gamma$. Then one extends to arbitrary $r$ by thinking of $\rho_g$ as a block matrix with blocks of size $|\Gamma| \times |\Gamma|$. 

Combining \eqref{eq:4.14} and \eqref{eq:4.15} we find:
\begin{eqnarray*} 
\log_p \ddet_{\Gamma} f & = & \tr_{\Gamma} \log f \\
& = & \frac{1}{|\Gamma|} \tr (\rho_{\log f})\\
& = & \frac{1}{|\Gamma|} \tr (\log \rho_f) \\
& = & \frac{1}{|\Gamma|} \log_p \ddet_{\Q_p} (\rho_f) \; .
\end{eqnarray*}
The last equation is proved by writing $\rho_f$ in triangular form in a suitable basis over $\oQ_p$ and observing that the eigenvalues of $\rho_f$ are $1$-units in $\oQ_p$. 
\end{proofof}

The next result is necessary to prove the relation of $\log_p \det_{\Gamma} f$ with $p$-adic entropies.

\begin{prop} \label{t:4.16}
Let $\Gamma$ be a residually finite countable discrete group and $\Gamma_n \to e$ a sequence as in the introduction. For $f$ in $1 + pM_r (c_0 (\Gamma , \Z_p))$ consider its image $f^{(n)}$ in $1 + pM_r (\Z_p \Gamma^{(n)})$ where $\Gamma^{(n)}$ is the finite group $\Gamma^{(n)} = \Gamma / \Gamma_n$. Then we have
\[
\log_p \ddet_{\Gamma} f = \lim_{n\to \infty} \log_p \ddet_{\Gamma^{(n)}} f^{(n)} \quad \mbox{in} \quad \Z_p \; .
\]
\end{prop}

\begin{proof}
The algebra map $M_r (c_0 (\Gamma)) \to M_r (c_0 (\Gamma^{(n)}))$ sending $f$ to $f^{(n)}$ is continuous since we have $\|f^{(n)} \| \le \|f \|$. For $f$ in $1 + p M_r (c_0 (\Gamma , \Z_p))$ we therefore get:
\[
(\log f)^{(n)} = \log f^{(n)} \quad \mbox{in} \; M_r (c_0 (\Gamma^{(n)})) \; .
\]
The next claim for $g = \log f$ thus implies the assertion.
\end{proof}

\begin{claim} \label{t:4.17}
For $g$ in $M_r (c_0 (\Gamma))$ we have
\[
\tr_{\Gamma} (g) = \lim_{n\to \infty} \tr_{\Gamma^{(n)}} (g^{(n)}) \; .
\]
\end{claim}

\begin{proof}
We may assume that $r = 1$. Writing $g = \sum a_{\gamma} \gamma$ with $a_{\gamma} \in \Q_p , |a_{\gamma}|_p \to 0$ for $\gamma \to \infty$ we have
\begin{eqnarray*}
|\tr_{\Gamma} (g) - \tr_{\Gamma^{(n)}} (g^{(n)})|_p & = & \Big| a_e - \sum_{\overline{\gamma} = e} a_{\gamma} \Big|_p \\
& = & \Big| \sum_{\gamma \in \Gamma_n \setminus e} a_{\gamma} \Big|_p \\
& \le & \max_{\gamma \in \Gamma_n \setminus e} |a_{\gamma}|_p \; .
\end{eqnarray*}
For $\varepsilon > 0$ there is a finite subset $S = S_{\varepsilon}$ of $\Gamma$ such that $|a_{\gamma}|_p < \varepsilon$ for $\gamma \in \Gamma \setminus S$. Only $e$ is contained in infinitely many $\Gamma_n$'s. Hence there is some $n_0$ such that $(\Gamma_n \setminus e) \cap S = \emptyset$ i.e. $\Gamma_n \setminus e \subset \Gamma \setminus S$ for all $n \ge n_0$. It follows that for $n \ge n_0$ we have $|\tr_{\Gamma} (g) - \tr_{\Gamma^{(n)}} (g^{(n)})|_p \le \varepsilon$. 
\end{proof}

\begin{cor} \label{t:4.18}
Let $\Gamma$ be a residually finite countable discrete group and $f$ an element of $M_r (\Z \Gamma)$ which is a $1$-unit in $M_r (c_0 (\Gamma))$. Then the $p$-adic entropy $h_p (f)$ of the $\Gamma$-action on $X_f$ exists for all $\Gamma_n \to e$ and we have
\[
h_p (f) = \log_p \ddet_{\Gamma} f \quad \mbox{in} \; \Z_p \; .
\]
\end{cor}

\begin{proof}
By propositions \ref{t:3.1} and \ref{t:4.12} we have
\begin{eqnarray*}
\frac{1}{(\Gamma : \Gamma_n)} \log_p |\Fix_{\Gamma_n} (X_f)| & = & \frac{1}{(\Gamma : \Gamma_n)} \log_p \ddet_{\Q_p} (\rho_{f^{(n)}}) \\
& = & \log_p \ddet_{\Gamma^{(n)}} (f^{(n)}) \; .
\end{eqnarray*}
Hence the assertion follows from proposition \ref{t:4.16}.
\end{proof}

\section{A $p$-adic logarithmic Fuglede--Kadison determinant and its relation to $p$-adic entropy} \label{sec:5}

Having defined a homomorphism $\log_p \det_{\Gamma}$ on $1 + p c_0 (\Gamma , \Z_p)$ in \eqref{eq:4.8} one would like to use the exact sequence \eqref{eq:2.10} to extend it to $c_0 (\Gamma , \Z_p)^*$. However, for infinite groups $\Gamma$ the abelianization of the group $\F_p [\Gamma]^*$ divided by the image of $\Gamma$ is not known to be torsion in any generality  -- as far as I know. However, corresponding results are known for $K_1$ of $\F_p [\Gamma]$ and this determines our approach which even for $r = 1$ requires the preceeding considerations for matrix algebras.

For a unital ring $R$ recall the embedding $\GL_r (R) \hookrightarrow \GL_{r+1} (R)$ mapping $a$ to $\left( \begin{smallmatrix} a & 0 \\ 0 & 1 \end{smallmatrix} \right)$. Let $\GL_{\infty} (R)$ be the union of the $\GL_r (R)$'s. We will view elements of $\GL_{\infty} (R)$ as infinite matrices with $1$'s on the diagonal and only finitely many further nonzero entries. The subgroup $E_r (R) \subset \GL_r (R)$ of elementary matrices is the subgroup generated by matrices which have $1$'s on the diagonal and at most one further non-zero entry. Let $E_{\infty} (R)$ be their union and set $K_1 (R) = \GL_{\infty} (R) / E_{\infty} (R)$. It is known that we have $E_{\infty} (R) = (\GL_{\infty} (R) , \GL_{\infty} (R))$ and hence that $K_1 (R) = \GL_{\infty} (R)^{\abb}$ c.f. \cite{M} \S\,3. The Whitehead group over $\F_p$ of a discrete group $\Gamma$ is defined to be 
\[
Wh^{\F_p} (\Gamma) := K_1 (\F_p [\Gamma]) / \langle \Gamma \rangle \; .
\]
Here $\langle \Gamma \rangle$ is the image of $\Gamma$ under the canonical map $\F_p [\Gamma]^* \to K_1 (\F_p [\Gamma])$.

We can treat groups for which $Wh^{\F_p} (\Gamma)$ is torsion. According to \cite{FL} Theorem 1.1 this is the case for torsion-free elementary amenable groups $\Gamma$. Recently, in \cite{L} it has been shown for a larger class of groups that $Wh^{\F_p} (\Gamma)$ is torsion. Apart from the elementary amenable groups, this class comprises all word hyperbolic groups. It is closed under subgroups, finite products, colimits and suitable extensions.

\begin{theorem} \label{t:5.1}
Let $\Gamma$ be a countable discrete residually finite group such that $Wh^{\F_p} (\Gamma)$ is torsion. Then there is a unique homomorphism
\[
\log_p \ddet_{\Gamma} : K_1 (c_0 (\Gamma , \Z_p)) \longrightarrow \Q_p
\]
with the following properties:\\
{\bf a} For every $r \ge 1$ the composition
\[
1 + pM_r (c_0 (\Gamma , \Z_p)) \hookrightarrow \GL_r (c_0 (\Gamma , \Z_p)) \to K_1 (c_0 (\Gamma , \Z_p)) \xrightarrow{\log_p \ddet_{\Gamma}} \Q_p
\]
coincides with the map $\log_p \det_{\Gamma}$ introduced in \eqref{eq:4.8}.\\
{\bf b} On the image of $\Gamma$ in $K_1 (c_0 (\Gamma , \Z_p))$ the map $\log_p \det_{\Gamma}$ vanishes.
\end{theorem}

\begin{proof}
Set $A = c_0 (\Gamma , \Z_p)$ and $\oA = A / pA = \F_p [\Gamma]$. The reduction map $A \to \oA$ induces an exact sequence
\begin{equation} \label{eq:5.2}
0 \to \Gamma E_{\infty} (A) (1 + p M_{\infty} (A)) / \Gamma E_{\infty} (A) \to K_1 (A) / \langle \Gamma \rangle\to K_1 (\oA) / \langle \Gamma \rangle \; .
\end{equation}
Here $M_{\infty} (A)$ is the (non-unital) algebra of infinite matrices $(a_{ij})_{i,j \ge 1}$ with only finitely many non-zero entries. Note that $1 + p M_{\infty} (A)$ is a subgroup of $\GL_{\infty} (A)$ since $1 + pM_r (A)$ is a subgroup of $\GL_r (A)$. Moreover $\Gamma E_{\infty} (A)$ is a normal subgroup of $\GL_{\infty} (A)$. Hence the sequence \eqref{eq:5.2} becomes an exact sequence:
\begin{equation} \label{eq:5.3}
0 \to (1 + pM_{\infty} (A)) / \Gamma E_{\infty} (A) \cap (1 + pM_{\infty} (A)) \to K_1 (A) / \langle \Gamma \rangle \to K_1 (\oA) / \langle \Gamma \rangle \; .
\end{equation}
Since $\Q_p$ is uniquely divisible this implies the uniqueness assertion in the theorem for any group $\Gamma$ such that $Wh^{\F_p} (\Gamma) = K_1 (\oA) / \langle \Gamma \rangle$ is torsion. For the existence, we first note that the homomorphisms defined in \eqref{eq:4.8} induce a homomorphism
\[
\log_p \ddet_{\Gamma} : 1 + p M_{\infty} (A) \longrightarrow \Z_p \; .
\]
We have to show that $\log_p \ddet_{\Gamma} f = 0$ for every $f$ in $1 + p M_{\infty} (A)$ which also lies in $\Gamma E_{\infty} (A)$. Under our identification of $\GL_r (A)$ with a subgroup of $\GL_{\infty} (A)$ we find some $r \ge 1$ such that we have
\[
f = i (\gamma) e_1 \cdots e_N \quad \mbox{in} \; 1 + pM_r (A) \; .
\]
Here the $e_i$ are elementary $r \times r$-matrices and $i (\gamma) = \left( \begin{smallmatrix} \gamma & 0 \\ 0 & 1_{r-1} \end{smallmatrix} \right)$ for some $\gamma$ in $\Gamma$.

According to propositions \ref{t:4.12} and \ref{t:4.16} we have for any choice of sequence $\Gamma_n \to e$:
\[
\log_p \ddet_{\Gamma} f = \lim_{n\to \infty} \frac{1}{(\Gamma : \Gamma_n)} \log_p \ddet_{\Q_p} (\rho_{f^{(n)}}) \; .
\]
On the other hand:
\[
\ddet_{\Q_p} (\rho_{f^{(n)}}) = \ddet_{\Q_p} (\rho_{i (\gamma)^{(n)} e^{(n)}_1 \ldots e^{(n)}_N}) = \ddet_{\Q_p} (\rho_{i (\gamma)^{(n)}}) \prod_i \ddet_{\Q_p} (\rho_{e^{(n)}_i}) \; .
\]
Let $\eb$ be a basis of $\Q_p [\Gamma^{(n)}]$. In the basis $(\eb , \ldots , \eb)$ of $\Q_p [\Gamma^{(n)}]^r$ the endomorphism $\rho_{e^{(n)}_i}$ is given by a matrix of $|\Gamma^{(n)}| \times |\Gamma^{(n)}|$-blocks. The diagonal blocks are identity matrices. At most one of the other blocks is non-zero. In particular, the matrix is triangular and we have $\ddet_{\Q_p} (\rho_{e^{(n)}_i}) = 1$. In the same basis $\rho_{i (\gamma)^{(n)}}$ is a permutation matrix and hence $\ddet_{\Q_p} (\rho_{i (\gamma)^{(n)}}) = \pm 1$. It follows that we have $\log_p \ddet_{\Gamma} f = 0$ as we wanted to show.
\end{proof}

I think that theorem \ref{t:5.1} should also hold without the condition that $\Gamma$ is residually finite.

\begin{remark} \label{t:5.4}
\rm For $\Gamma = \Z^d$ and $f$ in $\GL_r (c_0 (\Gamma , \Z_p))$, writing $[f]$ for the class of $f$ in $K_1$, we have
\[
\log_p \ddet_{\Gamma} [f] = m_p (\det f) 
\]
extending equation \eqref{eq:4.10}.
\end{remark}

This follows from the uniqueness assertion in theorem \ref{t:5.1}. Namely, the map $[f] \mapsto m_p (\det f)$ defines a homomorphism on $K_1$ which according to equation \eqref{eq:4.10} satisfies condition {\bf a}. It satisfies condition {\bf b} as well, since $\log_p$ vanishes on roots of unity and hence we have $m_p (t^{\nu}) = 0$ for all $\nu$ in $\Z^d$, c.f. example \ref{t:2.12}.

\begin{defn} \label{t:5.5}
For any group $\Gamma$ as in the theorem we define the homomorphism $\log_p \det_{\Gamma}$ on $\GL_r (c_0 (\Gamma , \Z_p))$ to be the composition
\[
\log_p \ddet_{\Gamma} : \GL_r (c_0 (\Gamma , \Z_p)) \longrightarrow K_1 (c_0 (\Gamma , \Z_p)) \xrightarrow{\log_p \det_{\Gamma}} \Q_p \; .
\]
\end{defn}

If we unravel the definitions we get the following description of this map. Given a matrix $f$ in $\GL_r (c_0 (\Gamma , \Z_p))$ there are integers $N \ge 1$ and $s \ge r$ such that in $\GL_s (c_0 (\Gamma , \Z_p))$ we have $f^N = i (\gamma) \varepsilon g$ with $\varepsilon$ in $E_s (c_0 (\Gamma , \Z_p))$, $g$ in $1 + p M_s (c_0 (\Gamma , \Z_p))$ and $i (\gamma)$ the $s \times s$-matrix $\left( \begin{smallmatrix} \gamma & 0 \\ 0 & 1_{s-1} \end{smallmatrix} \right)$. Then we have
\begin{equation} \label{eq:5.3_neu}
\log_p \ddet_{\Gamma} f = \frac{1}{N} \log_p \ddet_{\Gamma} g = \frac{1}{N} \tr_{\Gamma} \log g \; .
\end{equation}

We can now prove the following extension of corollary \ref{t:4.18}.

\begin{theorem} \label{t:5.6}
Let $\Gamma$ be a residually finite countable discrete group such that $Wh^{\F_p} (\Gamma)$ is torsion. Let $f$ be an element of $M_r (\Z \Gamma) \cap \GL_r (c_0 (\Gamma , \Z_p))$. Then $h_p (f)$ exists for all $\Gamma_n \to e$ and we have
\[
h_p (f) = \log_p \ddet_{\Gamma} f \quad \mbox{in} \; \Q_p \; .
\]
\end{theorem}

\begin{proof}
Let us write $f^N = i (\gamma) \varepsilon g$ as above. Then by proposition \ref{t:3.1} we have
\begin{eqnarray*}
\log_p |\Fix_{\Gamma_n} (X_f)| & = & \log_p \ddet_{\Q_p} (\rho_{f^{(n)}}) \\
& = & \frac{1}{N} \log_p \ddet_{\Q_p} (\rho_{i (\gamma)^{(n)}}) + \frac{1}{N} \log_p \ddet_{\Q_p} (\rho_{\varepsilon^{(n)}}) \\
& & + \frac{1}{N} \log_p \ddet_{\Q_p} (\rho_{g^{(n)}}) \; .
\end{eqnarray*}
Note here that the composition
\[
M_s (c_0 (\Gamma)) \longrightarrow M_s (c_0 (\Gamma^{(n)})) \xrightarrow{\rho} \End_{\Q_p} (\Q_p \Gamma^{(n)})^s
\]
is a homomorphism of algebras.

As in the proof of theorem \ref{t:5.1} we see that the terms $\log_p \ddet_{\Q_p} (\rho_{i (\gamma)^{(n)}})$ and $\log_p \det_{\Q_p} (\rho_{\varepsilon^{(n)}})$ vanish. This gives
\begin{eqnarray*}
\frac{1}{(\Gamma : \Gamma_n)} \log_p |\Fix_{\Gamma_n} (X_f)| & = & \frac{1}{(\Gamma : \Gamma_n)} \frac{1}{N} \log_p \ddet_{\Q_p} (\rho_{g^{(n)}}) \\
 & = & \frac{1}{N} \log_p \ddet_{\Gamma^{(n)}} (g^{(n)}) \quad \mbox{by proposition \ref{t:4.12}.}
\end{eqnarray*}
Using proposition \ref{t:4.16} we get in the limit $n\to \infty$ that
\[
h_p (f) = \frac{1}{N} \log_p \ddet_{\Gamma} (g) \overset{\eqref{eq:5.3_neu}}{=} \log_p \det f \; .
\]
\end{proof}

For groups $\Gamma$ as in theorem \ref{t:5.6} whose group ring $\F_p \Gamma$ has no zero divisors it is possible to extend the definition of $\log_p \det_{\Gamma}$ from $c_0 (\Gamma , \Z_p)^*$ to $c_0 (\Gamma)^*$. Namely, by proposition \ref{t:2.11} we know that
\[
c_0 (\Gamma)^* = p^{\Z} c_0 (\Gamma , \Z_p)^* \quad \mbox{and} \quad p^{\Z} \cap c_0 (\Gamma , \Z_p)^* = 1 \; .
\]
Hence there is a unique homomorphism
\[
\log_p \ddet_{\Gamma} : c_0 (\Gamma)^* \longrightarrow \Q_p
\]
which agrees with $\log_p \ddet_{\Gamma}$ previously defined on $c_0 (\Gamma , \Z_p)^*$ in definition \ref{t:5.5} and satisfies
\[
\log_p \ddet_{\Gamma} (p) = 0 \; .
\]
Let $\Gamma$ be a torsion-free elementary amenable group. Then according to \cite{KLM} Theorem 1.4 the group ring $\F_p \Gamma$ has no zero divisors and according to \cite{FL} Theorem 1.1 the group $Wh^{\F_p} (\Gamma)$ is torsion. Hence $\log_p \det_{\Gamma}$ is defined on $c_0 (\Gamma)^*$ and this is the map used in theorem \ref{t:1.9}. 

\begin{proofof} {\bf theorem \ref{t:1.9}}
Writing $f$ in $\Z \Gamma \cap c_0 (\Gamma)^*$ as a product $f = p^{\nu} g$ with $g$ in $c_0 (\Gamma , \Z_p)^*$ it follows that $g \in \Z \Gamma$ and proposition \ref{t:3.1} shows that we have
\begin{eqnarray*}
\log_p |\Fix_{\Gamma_n} (X_f)| & = & \log_p \ddet_{\Q_p} (\rho_{f^{(n)}}) \\
& = & \log_p \ddet_{\Q_p} (\rho_{g^{(n)}}) \\
 & = & \log_p |\Fix_{\Gamma_n} (X_g)| \; .
\end{eqnarray*}
Note here that we have $\log_p (p) = 0$. It follows from theorem \ref{t:5.6} applied to $g$ that for all $\Gamma_n \to e$ we get:
\[
h_p (f) = h_p (g) = \log_p \ddet_{\Gamma}  g = \log_p \ddet_{\Gamma} f \; .
\]
\end{proofof}

For $\Gamma = \Z^d$ it follows from remark \ref{t:5.4} that for any $f$ in $c_0 (\Z^d)^* = \Q_p \langle t^{\pm 1}_1 , \ldots , t^{\pm 1}_d \rangle^*$ we have:
\[
\log_p \ddet_{\Gamma} f = m_p (f) \; .
\]
Hence theorem \ref{t:1.7} is a special case of theorem \ref{t:1.9}.

Concerning approximations of $\log_p \det_{\Gamma} f$ we note that proposition \ref{t:4.16} extends to more general cases.

\begin{prop} \label{t:5.7}
Let $\Gamma$ be a residually finite countable discrete group and $\Gamma_n \to e$ as in the introduction. For $f$ in $M_r (c_0 (\Gamma))$ let $f^{(n)}$ be its image in $M_r (\Q_p \Gamma^{(n)})$. Then the formula
\begin{equation} \label{eq:5.8}
\log_p \ddet_{\Gamma} f = \lim_{n\to \infty} \frac{1}{(\Gamma : \Gamma_n)} \log_p \ddet_{\Q_p} (\rho_{f^{(n)}})
\end{equation}
holds whenever $\log_p \det_{\Gamma} f$ is defined. These are the cases\\
{\bf a} where $f$ is in $1 + p M_r (c_0 (\Gamma , \Z_p))$\\
{\bf b} where $Wh^{\F_p} (\Gamma)$ is torsion and $f$ is in $\GL_r (c_0 (\Gamma , \Z_p))$\\
{\bf c} where $Wh^{\F_p} (\Gamma)$ is torsion, $\F_p \Gamma$ has no zero divisors and $f$ is in $c_0 (\Gamma)^*$. 
\end{prop}

\begin{proof}
The assertions follow from propositions \ref{t:4.12} and \ref{t:4.16} together with calculations as in the proofs of theorems \ref{t:5.1} and \ref{t:5.6}.
\end{proof}

We end the paper with some open questions: Is there a dynamical criterion for the existence of the limit defining $p$-adic entropy? Is there a notion of ``$p$-adic expansiveness'' for $\Gamma$-actions on compact spaces $X$ which for the systems $X_f$ with $f$ in $M_r (\Z \Gamma)$ translates into the condition that $f$ is invertible in $M_r (c_0 (\Gamma))$, (or in $M_r (c_0 (\Gamma , \Z_p))$)? In fact, I assume that $p$-adic entropy can only be defined for ``$p$-adically expansive'' systems, c.f. \cite{BD} Remark after proposition 1.3. Is there a direct proof that the limit in formula \eqref{eq:5.8} exists?

Finally, in \cite{BD} a second version of a $p$-adic Mahler measure was defined which involves both the $p$-adic and the archimedian valuations of $\Q$. Can this be obtained for the systems $X_f$ by doing something more involved with the fixed points than taking their cardinalities and forming the limit \eqref{eq:1.3}?


Mathematisches Institut,
Einsteinstr. 62,
48149 M\"unster, Germany
\end{document}